\documentclass[11pt]{article}
\usepackage{latexsym,amssymb,color}
\usepackage{amsmath}
\newtheorem{Theorem}{\bf Theorem}[section]
\newtheorem{Lemma}{\bf Lemma}[section]
\newtheorem{Proposition}{\bf Proposition}[section]
\newtheorem{Corollary}{\bf Corollary}[section]
\newtheorem{Remark}{\bf Remark}[section]
\newtheorem{Example}{\bf Example}[section]
\newtheorem{Definition}{\bf Definition}[section]

\newenvironment{theorem}{\begin{Theorem}$\!\!\!$}{\end{Theorem}}
\newenvironment{lemma}{\begin{Lemma}$\!\!\!$}{\end{Lemma}}

\newenvironment{corollary}{\begin{Corollary}$\!\!\!$}{\end{Corollary}}

\newenvironment{definition}{\begin{Definition}$\!\!\!$}{\end{Definition}}

\numberwithin{equation}{section}

\begin{document}
\title{The large diffusion limit 
for the heat equation\\ 
with a dynamical boundary condition}
\author{
Marek Fila\\
Department of Applied Mathematics and Statistics\\
Comenius University,\\
84248 Bratislava, Slovakia
\\
\\
Kazuhiro Ishige\\
Graduate School of Mathematical Sciences\\ 
The University of Tokyo\\
3-8-1 Komaba, Meguro-ku, Tokyo 153-8914, Japan
\\
\\
Tatsuki Kawakami\\
Department of Applied Mathematics and Informatics,\\ 
Ryukoku University\\
Seta Otsu 520-2194, Japan
}
\date{}
\maketitle

\begin{abstract}
We study the heat equation on a half-space with a linear dynamical boundary condition. 
Our main aim is to show that, if the diffusion coefficient tends to infinity, then the solutions 
converge (in a suitable sense) to solutions of the Laplace equation with the same dynamical
boundary condition.
\end{abstract}
\vspace{10pt}
Keywords: heat equation, dynamical boundary condition,
 large diffusion limit

\section{Introduction}
We consider the problem
\begin{equation}
\label{eq:1.1}
\left\{
\begin{array}{ll}
\displaystyle{\varepsilon \partial_t u_\varepsilon-\Delta u_\varepsilon=0}, & x\in{\mathbb R}^N_+,\,\,\,t>0,\vspace{5pt}\\
\displaystyle{\partial_tu_\varepsilon+\partial_\nu u_\varepsilon=0}, & x\in\partial{\mathbb R}^N_+,\,\,\, t>0,\vspace{5pt}\\
\displaystyle{u_\varepsilon(x,0)=\varphi(x)},\qquad & x\in{\mathbb R}^N_+,
\vspace{5pt}\\
\displaystyle{u_\varepsilon(x,0)=\varphi_b(x')},\qquad & x=(x',0)\in\partial{\mathbb R}^N_+,
\end{array}
\right.
\end{equation}
where $N\ge 2$, ${\mathbb R}^N_+:={\mathbb R}^{N-1}\times{\mathbb R}_+$, 
$\Delta$ is the $N$-dimensional Laplacian (in $x$), 
$\partial_t:=\partial/\partial t$, $\partial_\nu:=-\partial/\partial x_N$, $\varepsilon\in(0,1)$ and 
$\varphi$ and $\varphi_b$ are measurable functions 
in ${\mathbb R}^N_+$ and ${\mathbb R}^{N-1}$, respectively. 
Our main aim is to show that, as $\varepsilon\to 0$, it holds that $u_\varepsilon\to u$ (in a suitable sense),
where $u$ is the solution of
\begin{equation}
\label{eq:1.2}
\left\{
\begin{array}{ll}
\displaystyle{\Delta u=0}, & x\in{\mathbb R}^N_+,\,\,\,t>0,\vspace{5pt}\\
\displaystyle{\partial_tu+\partial_\nu u=0}, & x\in\partial{\mathbb R}^N_+,\,\,\, t>0,\vspace{5pt}\\
\displaystyle{u(x,0)=\varphi_b(x')},\qquad & x=(x',0)\in\partial{\mathbb R}^N_+.
\end{array}
\right.
\end{equation}
This convergence does not look unexpected, see \cite{AAIY}, but we are not aware of any previous result which would
support this natural conjecture. In particular, convergence of this type means that the influence of the initial function
$\varphi$ is lost in the limit, and we shall describe this phenomenon in more detail.

A result in a similar spirit was established in \cite{AAIY} for the eikonal equation with the same dynamical boundary condition
as in \eqref{eq:1.1}. More precisely, the following problem was considered in \cite{AAIY}:
$$
\left\{
\begin{array}{ll}
\displaystyle{\varepsilon \partial_t u_\varepsilon+ |\nabla_x u_\varepsilon|=1}, & x\in\Omega,\,\,\,t>0,\vspace{5pt}\\
\displaystyle{\partial_tu_\varepsilon+\partial_\nu u_\varepsilon=0}, & x\in\partial\Omega,\,\,\, t>0,\vspace{5pt}\\
\displaystyle{u_\varepsilon(x,0)=\varphi(x)},\qquad & x\in\overline\Omega,
\end{array}
\right.
$$
where $\varepsilon\in(0,1)$, $\Omega\subset {\mathbb R}^N$ is a bounded domain with $C^1$-boundary, 
and $\nu$ is the outer normal of $\partial\Omega$. It was shown in \cite{AAIY} that $u_\varepsilon\to u$
as $\varepsilon\to 0$, where $u$ is the solution of
\[
\left\{
\begin{array}{ll}
\displaystyle{|\nabla_x u|=1}, & x\in\Omega,\,\,\,t>0,\vspace{5pt}\\
\displaystyle{\partial_tu+\partial_\nu u=0}, & x\in\partial\Omega,\,\,\, t>0,
\end{array}
\right.
\]
with an appropriate initial condition.

In the context of diffusion, the boundary condition from
\eqref{eq:1.1} describes thermal contact with a perfect conductor or
diffusion of solute from a well-stirred fluid or vapour (see e.g. \cite{C}). Various aspects of analysis of parabolic
and elliptic equations with dynamical boundary conditions have been treated by many authors, see for example
\cite{AQR}--\cite{BPR}, \cite{EMR, E2, E3, FQ2, FV, GH, H, KI}, \cite{VV1}--\cite{V1} for the parabolic case 
and \cite{E1}, \cite{E4}--\cite{FQ1}, \cite{GM}, \cite{KI}--\cite{L}, \cite{V2, Y} for the elliptic one. Here we
demonstrate on the simplest possible linear example how are these two classes of problems linked.

Throughout this paper we often identify ${\mathbb R}^{N-1}$ with $\partial{\mathbb R}^N_+$. 
We introduce some notation. 
Let $\Gamma_D=\Gamma_D(x,y,t)$ be the Dirichlet heat kernel on ${\mathbb R}^N_+$, that is,
\begin{equation}
\label{eq:1.3}
\Gamma_D(x,y,t):=(4\pi t)^{-\frac{N}{2}}
\left[\exp\biggr(-\frac{|x-y|^2}{4t}\biggr)-\exp\biggr(-\frac{|x-y_*|^2}{4t}\biggr)\right]
\end{equation}
for $(x,y,t)\in\overline{{\mathbb R}^{N}_+}\times{\mathbb R}^{N}_+\times(0,\infty)$,
where $y_*=(y',-y_N)$ for $y=(y',y_N)\in{\mathbb R}^N_+$.
Define
\begin{equation}
\label{eq:1.4}
[S_1(t)\phi](x):=\int_{{\mathbb R}^N_+}\Gamma_D(x,y,t)\phi(y)\, dy,
\quad (x,t)\in\overline{{\mathbb R}^N_+}\times(0,\infty),
\end{equation}
for any measurable function $\phi$ in ${\mathbb R}^N_+$. 
For $x=(x',x_N)\in\overline{{\mathbb R}^N_+}$ and $t>0$, 
set
$$
P(x',x_N,t):=c_N(x_N+t)^{1-N}\left(1+\left|\frac{x'}{x_N+t}\right|^2\right)^{-\frac{N}{2}},
$$
where $c_N$ is a constant chosen so that
$$
\int_{{\mathbb R}^{N-1}}P(x',x_N,t)\,dx'=1\quad\mbox{for all $x_N\ge 0$ and $t>0$}. 
$$
Then $P=P(x',x_N,t)$ is the fundamental solution of 
the Laplace equation in ${\mathbb R}^N_+$ 
with the homogeneous dynamical boundary condition, that is, 
$P$ satisfies 
$$
\left\{
\begin{array}{ll}
\displaystyle{-\Delta P=0}, & x\in{\mathbb R}^N_+,\,\,\,t>0,\vspace{5pt}\\
\displaystyle{\partial_tP+\partial_\nu P=0}, & x\in\partial{\mathbb R}^N_+,\,\,\, t>0,\vspace{5pt}\\
\displaystyle{P(x,0)=\delta(x')},\qquad & x=(x',0)\in\partial{\mathbb R}^N_+,
\end{array}
\right.
$$
where $\delta=\delta(\cdot)$ is the Dirac delta function on $\partial{\mathbb R}^N_+={\mathbb R}^{N-1}$.  
Define
\begin{equation}
\label{eq:1.5}
[S_2(t)\psi](x):=\int_{{\mathbb R}^{N-1}}P(x'-y',x_N,t)\psi(y')\, dy',
\quad (x,t)\in\overline{{\mathbb R}^N_+}\times(0,\infty),
\end{equation}
for any measurable function $\psi$ in ${\mathbb R}^{N-1}$. 
\vspace{5pt}

We formulate the definition of a solution of \eqref{eq:1.1} 
by the use of the two integral kernels $\Gamma_D$ and $P$.  
For simplicity, let $\varphi_b=\varphi_b(x')$ and $g=g(x',t)$ be continuous functions 
in ${\mathbb R}^{N-1}$ and ${\mathbb R}^{N-1}\times(0,\infty)$, respectively, such that 
$\varphi_b(x')$ and $g(x',t)$ decay rapidly as $|x'|\to\infty$. 
Then the function
\begin{equation}  
\begin{split}
\label{eq:1.6}
 & w(x,t)=w(x',x_N,t)
 :=[S_2(t)\varphi_b](x)+\int_0^t[S_2(t-s)g(s)](x)\,ds
\end{split}
\end{equation}
can be defined for $x=(x',x_N)\in {\mathbb R}^N_+$ and $t>0$ 
and it is a classical solution of the Cauchy problem for the Laplace equation with a nonhomogeneous dynamical boundary condition
\begin{equation}
\label{eq:1.7}
\left\{
\begin{array}{ll}
\displaystyle{-\Delta w=0},
& x\in{\mathbb R}^N_+,\quad t>0,\vspace{5pt}\\
\displaystyle{\partial_tw+\partial_\nu w=g},
& x\in\partial{\mathbb R}^N_+,\quad t>0,\vspace{5pt}\\
\displaystyle{w(x,0)=\varphi_b(x')},\quad
&x=(x',0)\in\partial{\mathbb R}^N_+.
\end{array}
\right.
\end{equation}
It follows from \eqref{eq:1.6} that
\begin{equation}
\label{eq:1.8}
\begin{split}
\partial_tw(x,t):= & \,\,\int_{{\mathbb R}^{N-1}}\partial_tP(x'-y',x_N,t)\varphi_b(y')\,dy'
\\
&\,\,\,\,
+\int_{{\mathbb R}^{N-1}}P(x'-y',x_N,0)g(y',t)\,dy'\\
 & \,\,\,\,\,\,
+\int_0^t\int_{{\mathbb R}^{N-1}}\partial_tP(x'-y',x_N,t-s)g(y',s)\,dy'\,ds
\end{split}
\end{equation}
for $x=(x',x_N)\in{\mathbb R}^N_+$ and $t\in(0,T)$. 
Set 
\begin{equation}
\label{eq:1.9}
\Phi(x):=\varphi(x)-[S_2(0)\varphi_b](x). 
\end{equation}
Then the function  
$$
v_\varepsilon(x,t)
:=[S_1(\varepsilon^{-1}t)\Phi](x)
-\int_0^t[S_1(\varepsilon^{-1}(t-s))]\partial_t w(s)](x)\, ds
$$
satisfies 
\begin{equation}
\label{eq:1.10}
\left\{
\begin{array}{ll}
\displaystyle{\varepsilon \partial_tv_\varepsilon=\Delta v_\varepsilon -\varepsilon\partial_t w
},\quad
& x\in{\mathbb R}^N_+,\quad t>0,\vspace{5pt}\\
\displaystyle{v_\varepsilon=0},
& x\in\partial{\mathbb R}^N_+,\quad t>0,\vspace{5pt}\\
\displaystyle{v_\varepsilon(x,0)=\Phi(x)},
&x\in{\mathbb R}^N_+.
\end{array}
\right.
\end{equation}
Let $\partial_{x_N}:=\partial/\partial x_N$. 
If $g_\varepsilon(x',t):=\partial_{x_N} v_\varepsilon(x',0,t)$ for $x'\in{\mathbb R}^{N-1}$, $t>0$, and
$w_\varepsilon$ is defined as in \eqref{eq:1.6} with $g_\varepsilon$ instead of $g$,
then it follows from \eqref{eq:1.7}, \eqref{eq:1.8} and \eqref{eq:1.10} that $v_\varepsilon$ and $w_\varepsilon$
satisfy
\begin{equation}
\label{eq:1.11}
\left\{
\begin{array}{ll}
\displaystyle{\varepsilon \partial_t v_\varepsilon=\Delta v_\varepsilon-\varepsilon F_1[\varphi_b]-\varepsilon F_2[v_\varepsilon]},
\qquad & x\in{\mathbb R}^N_+,\,\,\,t>0,\vspace{5pt}\\
\displaystyle{\Delta w_\varepsilon=0}, & x\in{\mathbb R}^N_+,\,\,\,t>0,\vspace{5pt}\\
\displaystyle{v_\varepsilon=0},\quad
\displaystyle{\partial_tw_\varepsilon-\partial_{x_N} w_\varepsilon=\partial_{x_N}v_\varepsilon}, & x\in\partial{\mathbb R}^N_+,\,\,\, t>0,\vspace{5pt}\\
\displaystyle{v_\varepsilon(x,0)=\Phi(x)},\qquad & x\in{\mathbb R}^N_+,\vspace{5pt}\\
\displaystyle{w_\varepsilon(x,0)=\varphi_b(x')}, & x=(x',0)\in\partial{\mathbb R}^N_+,
\end{array}
\right.
\end{equation}
where
\begin{align}
\label{eq:1.12}
F_1[\varphi_b](x,t) &
:=\int_{{\mathbb R}^{N-1}}\partial_t P(x'-y',x_N,t)\varphi_b(y')\,dy',\\
\begin{split}
\label{eq:1.13}
F_2[v](x,t) &
:=\int_{{\mathbb R}^{N-1}}P(x'-y',x_N,0)\partial_{x_N}v(y',0,t)\,dy'\\
 &
\quad
+\int_0^t\int_{{\mathbb R}^{N-1}}\partial_t P(x'-y',x_N,t-s)\partial_{x_N}v(y',0,s)\,dy'\,ds. 
\end{split}
\end{align}
Furthermore, the function $u_\varepsilon:=v_\varepsilon+w_\varepsilon$ is a classical solution of \eqref{eq:1.1}. 
Motivated by this observation, 
we formulate the definition of a solution of \eqref{eq:1.1} via problem~\eqref{eq:1.11}. 
\begin{definition}
\label{Definition:1.1}
Let $\varphi$ and $\varphi_b$ be measurable functions 
in ${\mathbb R}^N_+$ and ${\mathbb R}^{N-1}$, respectively.
Let $0<T\le\infty$ and 
$$
v_\varepsilon,\,\,\partial_{x_N}v_\varepsilon,\,\, w_\varepsilon\in C(\overline{{\mathbb R}^N_+}\times(0,T)). 
$$
We call $(v_\varepsilon,w_\varepsilon)$ a solution of \eqref{eq:1.11} in ${\mathbb R}^N_+\times(0,T)$
if $v_\varepsilon$ and $w_\varepsilon$ satisfy
\begin{align*}
\begin{split}
v_\varepsilon(x,t)
 & =[S_1(\varepsilon^{-1}t)\Phi](x)
-\int_0^t[S_1(\varepsilon^{-1}(t-s))F_1[\varphi_b](s)](x)\,ds\\
 & \qquad\qquad\qquad\qquad\quad\,\,\,
-\int_0^t[S_1(\varepsilon^{-1}(t-s))F_2[v_\varepsilon](s)](x)\,ds,
 \end{split}
 \notag
\\
w_\varepsilon(x,t)
&
=[S_2(t)\varphi_b](x)+\int_0^t[S_2(t-s)\partial_{x_N}v_\varepsilon(s)](x)\, ds,
\end{align*}
for $x\in\overline{{\mathbb R}^N_+}$ and $t\in(0,T)$. 
In the case when $T=\infty$, we call $(v_\varepsilon,w_\varepsilon)$ a global-in-time solution of \eqref{eq:1.11}
and $u_\varepsilon$ a global-in-time solution of \eqref{eq:1.1}.
\end{definition} 

We are ready to state the main results of this paper.
For $1\le r\le\infty$, we write 
$|\cdot|_{L^r}:=\|\cdot\|_{L^r(\partial{\mathbb R}^N_+)}$ and $\|\cdot\|_{L^r}:=\|\cdot\|_{L^r({{\mathbb R}^N_+})}$ 
for simplicity.
\begin{theorem}
\label{Theorem:1.1}
Let $N\ge 2$, $\varepsilon\in(0,1)$, $\varphi\in L^\infty({\mathbb R}^N_+)$ and $\varphi_b\in L^\infty({\mathbb R}^{N-1})$. 
Then problem~\eqref{eq:1.11} possesses a unique global-in-time solution 
$(v_\varepsilon,w_\varepsilon)$ 
satisfying 
\begin{equation}
\label{eq:1.14}
\sup_{0<t<T}\,\left[\|v_\varepsilon(t)\|_{L^\infty}
+(\varepsilon^{-1}t)^{\frac{1}{2}}\|\partial_{x_N}v_\varepsilon(t)\|_{L^\infty}
+\|w_\varepsilon(t)\|_{L^\infty}\right]<\infty
\end{equation}
for any $T>0$. 
Furthermore, 
$v_\varepsilon$ and $w_\varepsilon$ are bounded and smooth in $\overline{{\mathbb R}^N_+}\times I$ 
for any bounded interval $I\subset(0,\infty)$ and have the following properties for any $\tau>0$:
\begin{itemize}
  \item[{\rm (a)}]
  There exists $C(\tau)>0$ such that 
  \begin{equation*}
  \begin{split}
   &
  \sup_{0<t<\tau}\,\left[\|v_\varepsilon(t)\|_{L^\infty}
  +(\varepsilon^{-1}t)^{\frac{1}{2}}\|\partial_{x_N}v_\varepsilon(t)\|_{L^\infty}
  +\|w_\varepsilon(t)\|_{L^\infty}\right]
  \\
  & \qquad\qquad\qquad\qquad\qquad\qquad\qquad
  \le C_\tau\bigg(\|\varphi\|_{L^\infty}+|\varphi_b|_{L^\infty}\bigg);
  \end{split}
  \end{equation*}
  \item[{\rm (b)}] 
  $\displaystyle{\lim_{\varepsilon\to0}\sup_{0<t<\tau}\, t^{\frac{1}{2}}\|v_\varepsilon(t)\|_{L^\infty({\mathbb R}^{N-1}\times(0,L))}=0}$ 
  for any $L>0$;
  \item[{\rm (c)}] 
$\displaystyle{\lim_{\varepsilon\to0}\sup_{0<t<\tau}\,\|w_\varepsilon(t)-S_2(t)\varphi_b\|_{L^\infty}=0}$.
\end{itemize}
\end{theorem}
As a corollary of Theorem~\ref{Theorem:1.1}, 
we see that the solution $u_\varepsilon=v_\varepsilon+w_\varepsilon$ of \eqref{eq:1.1} 
converges to the solution $S_2(t)\varphi_b$ of \eqref{eq:1.2}.
\begin{corollary}
\label{Corollary:1.1}
Assume the same conditions as in Theorem~{\rm\ref{Theorem:1.1}}. 
Let 
$(v_\varepsilon,w_\varepsilon)$ be the solution given in Theorem~{\rm\ref{Theorem:1.1}}. 
Then $u_\varepsilon=v_\varepsilon+w_\varepsilon$ is a classical solution of \eqref{eq:1.1} and it satisfies 
$$
\lim_{\varepsilon\to 0}\sup_{\tau_1<t<\tau_2}\,\|u_\varepsilon(t)-S_2(t)\varphi_b\|_{L^\infty(K)}=0
$$
for any compact set $K$ in $\overline{{\mathbb R}^N_+}$ and $0<\tau_1<\tau_2<\infty$. 
\end{corollary} 
We prepare some useful lemmata in Section~2 and then we give a proof of Theorem~\ref{Theorem:1.1} in Section~3.
\section{Preliminaries}
In this section
we prove several lemmata on $S_1(t)\phi$, $F_1[\varphi_b]$ and $F_2[v]$.  
In what follows, by the letter $C$
we denote generic positive constants (independent of $x$ and $t$)
and they may have different values also within the same line. 
\vspace{3pt}

We first recall the following properties of $S_1(t)\phi$ 
(see e.g. \cite{GGS}).
\begin{itemize}
\item[(${\rm G_1}$)]
For any $1\le q\le r\le\infty$,
$$
\|S_1(t)\phi\|_{L^r}\le C t^{-\frac{N}{2}(\frac{1}{q}-\frac{1}{r})}\|\phi\|_{L^q},\qquad t>0,
$$
for all $\phi\in L^q({\mathbb R}^N_+)$.
In particular, if $q=r$, then
\begin{equation}
\label{eq:2.1}
\sup_{t>0}\,\|S_1(t)\phi\|_{L^r}\le \|\phi\|_{L^r}.
\end{equation}
\item[(${\rm G_2}$)]
Let $\phi\in L^q({\mathbb R}^N_+)$ with $1\le q\le\infty$. 
Then, for any $T>0$, 
$S_1(t)\phi$ is bounded and smooth in $\overline{{\mathbb R}^N_+}\times(T,\infty)$. 
\end{itemize} 
Furthermore, we have:
\begin{lemma}
\label{Lemma:2.1}
Let $\phi\in L^\infty({\mathbb R}^N_+)$. 
Then
\begin{equation}
\label{eq:2.2}
\sup_{t>0}\,t^{\frac{1}{2}}\|\partial_{x_N}[S_1(t)\phi]\|_{L^\infty}\le \|\phi\|_{L^\infty}.
\end{equation}
Furthermore,
\begin{equation}
\label{eq:2.3}
\lim_{\varepsilon\to0}\sup_{t>0}\,\,t^{\frac{1}{2}}
\|S_1(\varepsilon^{-1}t)\phi\|_{L^\infty({\mathbb R}^{N-1}\times(0,L))}=0
\quad\mbox{for any $L>0$}.
\end{equation}
\end{lemma}
\noindent
{\bf Proof.}
It follows from \eqref{eq:1.3} that
\begin{equation}
\label{eq:2.4}
\begin{split}
 & K(x,y,t):=\partial_{x_N}\Gamma_D(x,y,t)\\
 & =\Gamma_{N-1}(x'-y',t)\times
 \\
 &\qquad
 \times
\left(-\frac{x_N-y_N}{2t}\Gamma_1(x_N-y_N,t)+\frac{x_N+y_N}{2t}\Gamma_1(x_N+y_N,t)\right)
\end{split}
\end{equation}
for $(x,y,t)\in\overline{{\mathbb R}^N_+}\times{\mathbb R}^N_+\times(0,\infty)$, 
where $\Gamma_d$ $(d=1,2,\dots)$ is the Gauss kernel in $\mathbb R^d$. 
Then
\begin{equation}
\label{eq:2.5}
\begin{split}
&
\int_{{\mathbb R}^N_+}|K(x,y,t)|\,dy
\\
& 
\le\int_0^\infty \left(\frac{|x_N-y_N|}{2t}\Gamma_1(x_N-y_N,t)+\frac{x_N+y_N}{2t}\Gamma_1(x_N+y_N,t)\right)\,dy_N\\
 & =(\pi t)^{-\frac{1}{2}}\int_0^\infty 2\eta e^{-\eta^2}\,d\eta=(\pi t)^{-\frac{1}{2}}
\end{split}
\end{equation}
for $x\in\overline{{\mathbb R}^N_+}$ and $t>0$.
By \eqref{eq:1.4} and \eqref{eq:2.5} we have 
$$
|\partial_{x_N}[S_1(t)\phi](x)|
\le\int_{{\mathbb R}^N_+}|K(x,y,t)||\phi(y)|\,dy
\le t^{-\frac{1}{2}}\|\phi\|_{L^\infty}
$$
for $x\in\overline{{\mathbb R}^N_+}$ and $t>0$.
This implies \eqref{eq:2.2}.

On the other hand, 
for any $L>0$, it follows from \eqref{eq:1.3} that 
\begin{equation}
\label{eq:2.6}
\begin{split} 
& \int_{{\mathbb R}^N_+}\Gamma_D(x,y,\varepsilon^{-1}t)\,dy\\
& 
=\int_0^\infty\bigg(\Gamma_1(x_N-y_N,\varepsilon^{-1}t)-\Gamma_1(x_N+y_N,\varepsilon^{-1}t)\bigg)\,dy_N\\
& =2(4\pi\varepsilon^{-1}t)^{-\frac{1}{2}}\int_0^{x_N}\exp\left({-\frac{\varepsilon\eta^2}{4t}}\right)\,d\eta
\le 2(4\pi\varepsilon^{-1}t)^{-\frac{1}{2}}L\le C(\varepsilon^{-1} t)^{-\frac{1}{2}}
\end{split}
\end{equation}
for $x\in {\mathbb R}^{N-1}\times(0,L)$, $t>0$ and $\varepsilon>0$.
For any $\phi\in L^\infty({\mathbb R}^N_+)$,
by \eqref{eq:1.4} and \eqref{eq:2.6} we have
$$
|[S_1(\varepsilon^{-1}t)\phi](x)|
\le
\int_{{\mathbb R}^N_+}\Gamma_D(x,y,\varepsilon^{-1}t)|\phi(y)|\,dy
\le C(\varepsilon^{-1} t)^{-\frac{1}{2}}\|\phi\|_{L^\infty}
$$
for $x\in {\mathbb R}^{N-1}\times(0,L)$, $t>0$ and $\varepsilon>0$.
This implies \eqref{eq:2.3},
and the proof of Lemma~\ref{Lemma:2.1} is complete.
$\Box$ 
\vspace{5pt}

Next we recall some properties of $S_2(t)\psi$. 
\begin{itemize}
  \item[(${\rm P}_1$)] 
  Let $\psi\in L^r({\mathbb R}^{N-1})$ for some $r\in[1,\infty]$ 
  and $t$, $t'>0$. Then
  \begin{align*}
    & [S_2(t)\psi](x',x_N)=[S_2(t+x_N)\psi](x',0),\vspace{3pt}
    \\
    & [S_2(t+t')\psi](x)=[S_2(t)(S_2(t')\psi)](x),
  \end{align*}
  for $x=(x',x_N)\in\overline{{\mathbb R}^N_+}$. 
  Furthermore, 
  $$
  \lim_{t\to 0}|S_2(t)\psi-\psi|_r=0\quad\mbox{if $1\le r<\infty$};
  $$
  \item[(${\rm P}_2$)] 
  For any $1\le q\le r\le\infty$,  
 $$
  |S_2(t)\psi|_{L^r}\le Ct^{-(N-1)(\frac{1}{q}-\frac{1}{r})}|\psi|_{L^q},\qquad t>0, 
 $$ 
  for all $\psi\in L^q({\mathbb R}^{N-1})$. 
  In particular, if $q=r$, then 
  \begin{equation}
  \label{eq:2.7}
  \sup_{t>0}\,|S_2(t)\psi|_{L^r}\le |\psi|_{L^r}.
  \end{equation}
 \end{itemize}
Properties~(${\rm P}_1$) and (${\rm P}_2$) easily follow from \eqref{eq:1.5} (see e.g. \cite{FIK02}) and imply that 
\begin{align}
\label{eq:2.8}
\sup_{t>0}\,\|S_2(t)\psi\|_{L^\infty}\le |\psi|_{L^\infty}
\end{align}
for all $\psi\in L^\infty({\mathbb R}^{N-1})$. 
Furthermore,
by a similar argument as in the proof of property~($G_2$) 
we have: 
\begin{itemize}
  \item[(${\rm P}_3$)] 
  Let $\psi\in L^q({\mathbb R}^{N-1})$ with $1\le q\le\infty$. 
  Then, for any $T>0$, $S_2(t)\psi$ is bounded and smooth in $\overline{{\mathbb R}^N_+}\times(T,\infty)$. 
\end{itemize}
\begin{lemma}
\label{Lemma:2.2}
Let $\psi\in L^\infty({\mathbb R}^{N-1})$. 
Set 
\begin{equation}
\label{eq:2.9}
D_\varepsilon[\psi](x,t)
:=\int_0^t[S_1(\varepsilon^{-1}(t-s))F_1[\psi](s)](x)\,ds
\end{equation}
for $x\in\overline{{\mathbb R}^N_+}$, $t>0$ and $\varepsilon>0$. 
Then $D_\varepsilon[\psi]$ and $\partial_{x_N}D_\varepsilon[\psi]$ are bounded and smooth 
in $\overline{{\mathbb R}^N_+}\times(T,\infty)$ for any $T>0$. 
Furthermore, there exists $C>0$ such that 
\begin{equation}
\label{eq:2.10}
\|D_\varepsilon[\psi](t)\|_{L^\infty}\le 
Ct^{\frac{1}{4}}(\varepsilon^{\frac{1}{2}}+t^{\frac{3}{4}})|\psi|_{L^\infty}
\end{equation}
for $t>0$ and $\varepsilon>0$.
Moreover, 
\begin{equation}
\label{eq:2.11}
\lim_{\varepsilon\to0}\sup_{t\in(0,T_1)}\,
\|D_\varepsilon[\psi](t)\|_{L^\infty({\mathbb R}^{N-1}\times(0,L))}=0
\end{equation}
for $T_1>0$ and $L>0$. 
\end{lemma}
{\bf Proof.}
We prove \eqref{eq:2.10} first. 
Since 
$$
\partial_tP(x',x_N,t)
=\frac{1}{x_N+t}\frac{|x'|^2-(N-1)(x_N+t)^2}{|x'|^2+(x_N+t)^2}P(x',x_N,t), 
$$
it follows that 
\begin{equation}
\label{eq:2.12}
|\partial_tP(x',x_N,t)|
\le C(x_N+t)^{-1}P(x',x_N,t). 
\end{equation}
By \eqref{eq:1.12}, \eqref{eq:2.7} and \eqref{eq:2.12} we have
\begin{equation}
\label{eq:2.13}
\begin{split}
\|F_1[\psi](\cdot,y_N,s)\|_{L^\infty({\mathbb R}^{N-1})}
&
\le C(y_N+s)^{-1}\|S_2(s+y_N)\psi\|_{L^\infty({\mathbb R}^{N-1})}
\\
&
\le C(y_N+s)^{-1}|\psi|_{L^\infty}
\end{split}
\end{equation}
for $y_N\in[0,\infty)$ and $s>0$. 
Since
\begin{equation}
\label{eq:2.14}
(y_N+s)^{-1}\le
\left\{
\begin{array}{ll}
y_N^{-\frac{3}{4}}s^{-\frac{1}{4}} 
&
\quad\mbox{for}\quad 0\le y_N\le 1,
\vspace{5pt}
\\
1
&
\quad\mbox{for}\quad y_N>1,
\end{array}
\right.
\end{equation}
by \eqref{eq:1.3}, \eqref{eq:2.9} and \eqref{eq:2.13} 
we see that
\begin{equation}
\label{eq:2.15}
\begin{split}
 & |D_\varepsilon[\psi](x,t)|
 \le\int_0^t\int_{{\mathbb R}^N_+}\Gamma_D(x,y,\tau_\varepsilon)|F_1[\psi](y,s)|\,dy\,ds\\
 &
 \le 
 C \int_0^t\int_0^\infty\Gamma_1(x_N-y_N,\tau_\varepsilon)\|F_1[\psi](\cdot,y_N,s)\|_{L^\infty({\mathbb R}^{N-1})}\,dy_N\,ds\\
 &
 \le C|\psi|_{L^\infty}\int_0^t\int_0^\infty\tau_\varepsilon^{-\frac{1}{2}}
 \exp\left(-\frac{(x_N-y_N)^2}{4\tau_\varepsilon}\right)
 \left(y_N+s\right)^{-1}\,dy_N\,ds\\
 &
 \le C|\psi|_{L^\infty}\left\{
 \int_0^t\tau_\varepsilon^{-\frac{1}{2}}s^{-\frac{1}{4}}\int_0^1y_N^{-\frac{3}{4}}\,dy_N\,ds+\int_0^t\,ds\right\}\\
 &
 \le C|\psi|_{L^\infty}\left\{\varepsilon^{\frac{1}{2}}
 \int_0^t(t-s)^{-\frac{1}{2}}s^{-\frac{1}{4}}\,ds+t\right\}\\
 &
 \le C|\psi|_{L^\infty}(\varepsilon^{\frac{1}{2}}t^{\frac{1}{4}}+t)
 = C|\psi|_{L^\infty}t^{\frac{1}{4}}(\varepsilon^{\frac{1}{2}}+t^{\frac{3}{4}})
\end{split}
\end{equation}
for $x\in\overline{{\mathbb R}^N_+}$, $t>0$ and $\varepsilon>0$,
where $\tau_\varepsilon:=\varepsilon^{-1}(t-s)$.
Here $\Gamma_1$ is the one-dimensional Gauss kernel.
This implies \eqref{eq:2.10}. 

We prove \eqref{eq:2.11}. 
Let $L>0$. 
Similarly to \eqref{eq:2.15}, 
we obtain
\begin{equation}
\label{eq:2.16}
\begin{split}
 & |D_\varepsilon[\psi](x,t)|
 \le\int_0^t\int_{{\mathbb R}^N_+}\Gamma_D(x,y,\tau_\varepsilon)|F_1[\psi](y,s)|\,dy\,ds\\
 & \le C
 \int_0^t\int_0^\infty\bigg(\Gamma_1(x_N-y_N,\tau_\varepsilon)-\Gamma_1(x_N+y_N,\varepsilon^{-1}(t-s))\bigg)\\
 & \qquad\qquad\qquad\qquad\qquad\qquad
 \times
 \|F_1[\psi](\cdot,y_N,s)\|_{L^\infty({\mathbb R}^{N-1})}\,dy_N\,ds\\
 & \le C|\psi|_{L^\infty}
 \int_0^t\tau_\varepsilon^{-\frac{1}{2}}
 \int_0^\infty
 (y_N+s)^{-1}
 \\
 &\hspace{1cm}
 \times
 \left[\exp\left({-\frac{(x_N-y_N)^2}{4\tau_\varepsilon}}\right)-\exp\left({-\frac{(x_N+y_N)^2}{4\tau_\varepsilon}}\right)\right]
\,dy_N\,ds
 \end{split}
\end{equation}
for $x\in\overline{{\mathbb R}^N_+}$, $t>0$ and $\varepsilon>0$.
This together with \eqref{eq:2.6}, \eqref{eq:2.13} and \eqref{eq:2.14} 
implies that
\begin{equation*}
\begin{split}
 & |D_\varepsilon[\psi](x,t)|
 \\
 & \le C|\psi|_{L^\infty}\bigg\{\int_0^t\tau_\varepsilon^{-\frac{1}{2}}s^{-\frac{1}{4}}\int_0^1y_N^{-\frac{3}{4}}\,dy_N\,ds\\
 & \,\,\,
 +\int_0^t\tau_\varepsilon^{-\frac{1}{2}}
 \int_1^\infty
 \left[\exp\left({-\frac{(x_N-y_N)^2}{4\tau_\varepsilon}}\right)-\exp\left({-\frac{(x_N+y_N)^2}{4\tau_\varepsilon}}\right)\right]
 \, dy_N \,ds\bigg\}\\
 & \le C|\psi|_{L^\infty}\left\{\int_0^t\tau_\varepsilon^{-\frac{1}{2}}s^{-\frac{1}{4}}\,ds
 +\int_0^t\tau_\varepsilon^{-\frac{1}{2}}\,ds\right\}\\
 & \le C|\psi|_{L^\infty}\varepsilon^{\frac{1}{2}}\left\{\int_0^t(t-s)^{-\frac{1}{2}}s^{-\frac{1}{4}}\,ds
 +\int_0^t(t-s)^{-\frac{1}{2}}\,ds\right\}\\
 & \le C|\psi|_{L^\infty}(\varepsilon^{\frac{1}{2}}t^{\frac{1}{4}}+(\varepsilon t)^\frac{1}{2})
\end{split}
\end{equation*}
for $x\in {\mathbb R}^{N-1}\times(0,L)$, $t>0$ and $\varepsilon>0$.
Thus \eqref{eq:2.11} holds.

On the other hand, 
it follows from the semigroup property of $S_1(t)$ that 
\begin{equation*}
\begin{split}
 & D_\varepsilon[\psi](x,t)
 =\int_0^t[S_1(\varepsilon^{-1}(t-s))F_1[\psi](s)](x)\,ds\\
 & =S_1(\varepsilon^{-1}(t-T/2))D_\varepsilon[\psi](x,T/2)
+\int_{T/2}^t[S_1(\varepsilon^{-1}(t-s))F_1[\psi](s)](x)\,ds
\end{split}
\end{equation*}
for $x\in\overline{{\mathbb R}^N_+}$ and $0<T<t<\infty$. 
Then, by \eqref{eq:2.9} and (${\rm G_2}$) 
we see that 
$$
S_1(\varepsilon^{-1}(t-T/2))D_\varepsilon[\psi](x,T/2)
$$
is bounded and smooth in $\overline{{\mathbb R}^N_+}\times(T,\infty)$. 
Furthermore, by \eqref{eq:2.13} 
we apply the same argument as in \cite[Section~3, Chapter 1]{F} 
to see that 
$$
\int_{T/2}^t[S_1(\varepsilon^{-1}(t-s))F_1[\psi](s)](x)\,ds
$$
is also bounded and smooth in $\overline{{\mathbb R}^N_+}\times(T,\infty)$. 
(See also \cite[Proposition~5.2]{FIK01} and \cite[Lemma~2.1]{IKK01}.) 
Therefore we deduce that 
$D_\varepsilon[\psi]$ and $\partial_{x_N}D_\varepsilon[\psi]$ are bounded and smooth 
in $\overline{{\mathbb R}^N_+}\times(T,\infty)$.
Thus Lemma~\ref{Lemma:2.2} follows.  
$\Box$\vspace{5pt}
\begin{lemma}
\label{Lemma:2.3}
Let $0\le \alpha<1$. Then there exists $C>0$ such that 
\begin{equation}
\label{eq:2.17}
\int_0^\infty
\frac{|x\pm y|}{t}\Gamma_1(x\pm y,t)y^{-\alpha}\,dy
\le Ct^{-\frac{\alpha+1}{2}}
\end{equation}
for $x\ge 0$ and $t>0$.
Here $\Gamma_1$ is the one-dimensional Gauss kernel.
\end{lemma}
{\bf Proof.}
Let $x\ge 0$ and $t>0$. 
It follows that 
\begin{equation*}
\begin{split}
 & \int_0^\infty
\frac{|x-y|}{t}\Gamma_1(x-y,t)y^{-\alpha}\,dy\\
 & =(4\pi t)^{-\frac{1}{2}}\biggr(\int_0^{x_2}+\int_{x/2}^\infty\biggr)
\frac{|x-y|}{t}\exp\left(-\frac{|x-y|^2}{4t}\right)y^{-\alpha}\,dy. 
\end{split}
\end{equation*}
Since $y^{-1}\le |x-y|^{-1}$ for $0\le x\le 2y$,
we have
\begin{equation}
\label{eq:2.18}
\begin{split}
 & \int_0^\infty
\frac{|x-y|}{t}\Gamma_1(x-y,t)y^{-\alpha}\,dy\\
 & \le Ct^{-\frac{1}{2}}\frac{x}{t}
 \exp\left(-\frac{x^2}{16t}\right)
 \int_0^{x/2}y^{-\alpha}\,dy_N\\
 & \qquad\quad
 +Ct^{-\frac{3}{2}+\frac{1-\alpha}{2}}\int_{x/2}^\infty
 \left(\frac{|x-y|}{t^{1/2}}\right)^{1-\alpha}
 \exp\left(-\frac{|x-y|^2}{4t}\right)\,dy\\
 & \le Ct^{-\frac{3}{2}}x^{2-\alpha}\exp\left(-\frac{x^2}{16t}\right)
 +Ct^{-1+\frac{1-\alpha}{2}}\le Ct^{-\frac{\alpha+1}{2}}.
\end{split}
\end{equation}
Since $y^{-1}\le 2(x+y)^{-1}$ for $0\le x\le y$,
similarly to \eqref{eq:2.18}, we obtain
$$
\int_0^\infty
\frac{x+y}{t}\Gamma_1(x+y,t)y^{-\alpha}\,dy
\le Ct^{-\frac{\alpha+1}{2}}.
$$
Thus \eqref{eq:2.17} holds and Lemma~\ref{Lemma:2.3} follows.
$\Box$\vspace{5pt}

\begin{lemma}
\label{Lemma:2.4}
Let $\psi\in L^\infty({\mathbb R}^{N-1})$. 
Then there exists $C>0$ such that 
\begin{equation}
\label{eq:2.19}
\|\partial_{x_N}D_\varepsilon[\psi](t)\|_{L^\infty}
\le 
C\varepsilon^{\frac{3}{4}}t^{-\frac{1}{4}}|\psi|_{L^\infty}
\end{equation}
for $t>0$ and $\varepsilon>0$.
\end{lemma}
{\bf Proof.}
By \eqref{eq:2.4}, \eqref{eq:2.9} and \eqref{eq:2.13}
we see that
\begin{equation}
\label{eq:2.20}
\begin{split}
& |\partial_{x_N} D_\varepsilon[\psi](x,t)|\\
& \le\int_0^t\int_{{\mathbb R}^N_+}|K(x,y,\tau_\varepsilon)||F_1[\psi](y,s)|\,dy\,ds\\
& \le C
\int_0^t\int_0^\infty \tilde{K}(x_N,y_N,\tau_\varepsilon)
\|F_1[\psi](\cdot,y_N,s)\|_{L^\infty({\mathbb R}^{N-1})}\, dy_N\,ds\\
&
\le C|\psi|_{L^\infty}
\int_0^ts^{-\frac{1}{2}}\int_0^\infty
\tilde{K}(x_N,y_N,\tau_\varepsilon)
y_N^{-\frac{1}{2}}\, dy_N\,ds
\end{split}
\end{equation}
for $x\in\overline{{\mathbb R}^N_+}$, $t>0$ and $\varepsilon>0$,
where $\tau_\varepsilon:=\varepsilon^{-1}(t-s)$ and 
\begin{equation}
\label{eq:2.21}
\tilde{K}(x_N,y_N,t)=\frac{|x_N-y_N|}{t}\Gamma_1(x_N-y_N,t)
+\frac{x_N+y_N}{t}\Gamma_1(x_N+y_N,t)
\end{equation}
for $x_N\ge0$, $y_N>0$ and $t>0$.
By \eqref{eq:2.17} with $\alpha=1/2$ and \eqref{eq:2.20}
we deduce that
\begin{equation*}
\begin{split}
|\partial_{x_N} D_\varepsilon[\psi](x,t)|
&
\le C|\psi|_{L^\infty}
\int_0^ts^{-\frac{1}{2}}\tau_\varepsilon^{-\frac{3}{4}}\,ds
\\
&
= C|\psi|_{L^\infty}
\int_0^ts^{-\frac{1}{2}}(\varepsilon^{-1}(t-s))^{-\frac{3}{4}}\,ds
\le C\varepsilon^{\frac{3}{4}}t^{-\frac{1}{4}}|\psi|_{L^\infty}
\end{split}
\end{equation*}
for $x\in\overline{{\mathbb R}^N_+}$, $t>0$ and $\varepsilon>0$.
Thus \eqref{eq:2.19} follows.
$\Box$
\section{Proof of Theorem~\ref{Theorem:1.1}}
We introduce some notation. 
Let $T>0$ and $\varepsilon\in(0,1)$.  
Set 
$$
X_T:=\bigg\{v, \partial_{x_N}v\in C(\overline{{\mathbb R}^N_+}\times(0,T))\,:\,\|v\|_{X_T}<\infty\bigg\},
\quad
\|v\|_{X_T}:=\sup_{0<t<T}E_\varepsilon[v](t),
$$
where 
$$
E_\varepsilon[v](t):=\|v(t)\|_{L^\infty}
+(\varepsilon^{-1}t)^{\frac{1}{2}}\|\partial_{x_N}v(t)\|_{L^\infty}.
$$
Then $X_T$ is a Banach space equipped with the norm $\|\cdot\|_{X_T}$. 
For the proof of Theorem~\ref{Theorem:1.1} 
we apply the Banach contraction mapping principle in $X_T$ 
to find a fixed point of 
\begin{equation}
\label{eq:3.1}
Q_\varepsilon[v](t)
:=S_1(\varepsilon^{-1}t)\Phi
-D_\varepsilon[\varphi_b](t)-\int_0^tS_1(\varepsilon^{-1}(t-s))F_2[v](s)\,ds,
\end{equation}
where $\Phi$, $F_2[v]$ and $D_\varepsilon[\varphi_b]$ are as in \eqref{eq:1.9}, \eqref{eq:1.13} and \eqref{eq:2.9}, respectively. 
\begin{lemma}
\label{Lemma:3.1}
There exists $C>0$ such that 
\begin{equation}
\label{eq:3.2}
F_2[v](x,t)\le C\varepsilon^\frac{1}{2}\left(t^{-\frac{1}{2}}+h(x_N,t)\right)\|v\|_{X_T}
\end{equation}
for $x\in{\mathbb R}^N_+$, $0<t<T$, $\varepsilon\in(0,1)$ and $v\in X_T$. 
Here
$$
h(x_N,t):=x_N^{-\frac{3}{4}}t^{\frac{1}{4}}\quad\mbox{if}\quad 0<x_N\le1
\quad\mbox{and}\quad
h(x_N,t):=x_N^{-\frac{1}{2}}\quad\mbox{if}\quad x_N>1.
$$
\end{lemma}
\noindent
{\bf Proof.}
Let $T>0$, $\varepsilon\in(0,1)$ and $v\in X_T$. 
It follows from \eqref{eq:1.13} that 
\begin{equation}
\label{eq:3.3}
F_2[v](x,t)=F_2'[v](x,t)+F_2''[v](x,t)
\end{equation}
for $x\in{\mathbb R}^N_+$ and $t>0$, 
where 
\begin{align*}
 & F_2'[v](x,t):=\int_{{\mathbb R}^{N-1}}P(x'-y',x_N,0)\partial_{x_N}v(y',0,t)\, dy',\\
 & F_2''[v](x,t):=\int_0^t\int_{{\mathbb R}^{N-1}}\partial_t P(x'-y',x_N,t-s)\partial_{x_N}v(y',0,s)\, dy'\,ds.
\end{align*}
Since $v\in X_T$, 
by \eqref{eq:1.5} and \eqref{eq:2.7} we see that 
\begin{equation}
\label{eq:3.4}
\begin{split}
|F_2'[v](x,t)|
&
\le \|S_2(x_N)\partial_{x_N}v(\cdot,0,t)\|_{L^\infty({\mathbb R}^{N-1})}
\\
&
\le |\partial_{x_N}v(t)|_{L^\infty}
\le \varepsilon^{\frac{1}{2}}t^{-\frac{1}{2}}\|v\|_{X_T}
\end{split}
\end{equation}
for $x\in{\mathbb R}^N_+$ and $t>0$. 
On the other hand, 
it follows from \eqref{eq:2.12} that
$$
|\partial_t P(x',x_N,t)|
\le\left\{
\begin{array}{ll}
CP(x',x_N,t)x_N^{-\frac{3}{4}}t^{-\frac{1}{4}} \quad & \mbox{if}\quad x_N\le1,
\vspace{5pt}\\
CP(x',x_N,t)x_N^{-\frac{1}{2}}t^{-\frac{1}{2}} & \mbox{if}\quad x_N>1,
\end{array}
\right.
$$
for $x\in{\mathbb R}^N_{+}$ and $t>0$. 
Then we obtain
\begin{equation}
\label{eq:3.5}
\begin{split}
&
|F_2''[v](x,t)|
\\
&
\le Cx_N^{-\frac{3}{4}}\int_0^t(t-s)^{-\frac{1}{4}}
\int_{{\mathbb R}^{N-1}}P(x'-y',x_N,t-s)|\partial_{x_N}v(y',0,s)|\,dy'\,ds
\\
&
\le C\varepsilon^{\frac{1}{2}}x_N^{-\frac{3}{4}}\bigg(\sup_{0<s<t}\,(\varepsilon^{-1}s)^{\frac{1}{2}}|\partial_{x_N}v(s)|_{L^\infty}\bigg)
\int_0^t(t-s)^{-\frac{1}{4}}s^{-\frac{1}{2}}\,ds
\\
&
\le C\varepsilon^{\frac{1}{2}}x_N^{-\frac{3}{4}}t^{\frac{1}{4}}\|v\|_{X_T}
=C\varepsilon^{\frac{1}{2}}h(x_N,t)\|v\|_{X_T}
\end{split}
\end{equation}
for $x'\in{\mathbb R}^{N-1}$, $0<x_N\le 1$ and $0<t<T$. 
Similarly, we deduce that
\begin{equation}
\label{eq:3.6}
\begin{split}
&
|F_2''[v](x,t)|
\\
&
\le C\int_0^t(t-s)^{-\frac{1}{2}}\int_{{\mathbb R}^N_+}P(x'-y',x_N,t-s)|\partial_{x_N}v(y',0,s)|\,dy'\,ds\\
 &
\le C\varepsilon^{\frac{1}{2}}x_N^{-\frac{1}{2}}
\bigg(\sup_{0<s<t}\,(\varepsilon^{-1}s)^{\frac{1}{2}}|\partial_{x_N}v(s)|_{L^\infty}\bigg)
\int_0^t(t-s)^{-\frac{1}{2}}s^{-\frac{1}{2}}\,ds\\
 & 
\le C\varepsilon^{\frac{1}{2}}x_N^{-\frac{1}{2}}\|v\|_{X_T}
=C\varepsilon^{\frac{1}{2}}h(x_N,t)\|v\|_{X_T}
\end{split}
\end{equation}
for $x'\in{\mathbb R}^{N-1}$, $x_N>1$ and $0<t<T$.
Therefore, by \eqref{eq:3.3}, \eqref{eq:3.4}, \eqref{eq:3.5} and \eqref{eq:3.6} we obtain \eqref{eq:3.2}.
Thus Lemma~\ref{Lemma:3.1} follows. 
$\Box$
\begin{lemma}
\label{Lemma:3.2}
For any $v\in X_T$ and $\varepsilon\in(0,1)$,
set
\begin{equation}
\label{eq:3.7}
\tilde{D}_\varepsilon[v](t):=
\int_0^tS_1(\varepsilon^{-1}(t-s))F_2[v](s)\,ds.
\end{equation}
Then there exists $T_*=T_*(N)>0$ such that 
\begin{equation}
\label{eq:3.8}
\|\tilde{D}_\varepsilon[v]\|_{X_{T_*}}\le\frac{1}{4}\|v\|_{X_{T_*}}
\end{equation}
for $v\in X_{T_*}$ and $\varepsilon\in(0,1)$.
Furthermore,
$\tilde{D}_\varepsilon[v]$ and $\partial_{x_N}\tilde{D}_\varepsilon[v]$ are bounded and smooth in $\overline{{\mathbb R}^N_+}\times(\tau,T_*)$ for any $0<\tau<T_*$.
\end{lemma}
\noindent
{\bf Proof.}
Let $T>0$.
By \eqref{eq:2.9} and \eqref{eq:3.7} we see that 
$\tilde{D}_\varepsilon$ is defined analogously as  $D_\varepsilon$ with $F_1$ replaced by $F_2$. 
Then it follows from \eqref{eq:2.15} and \eqref{eq:3.2} that
\begin{equation*}
\begin{split}
& |\tilde{D}_\varepsilon[v](x,t)|\\
& \le C\int_0^t\int_0^\infty\Gamma_1(x_N-y_N,\tau_\varepsilon)\|F_2[v](\cdot,y_N,s)\|_{L^\infty({\mathbb R}^{N-1})}\,dy_N\,ds\\
& \le C\varepsilon^{\frac{1}{2}}\|v\|_{X_T}\int_0^t\int_0^\infty\tau_\varepsilon^{-\frac{1}{2}}\exp\left(-\frac{(x_N-y_N)^2}{4\tau_\varepsilon}\right)\left(s^{-\frac{1}{2}}+h(y_N,s)\right)\,dy_N\,ds\\
& \le C\varepsilon^{\frac{1}{2}}\|v\|_{X_T}\left\{\int_0^ts^{-\frac{1}{2}}\,ds
+\int_0^t\tau_\varepsilon^{-\frac{1}{2}}s^{\frac{1}{4}}\int_0^1y_N^{-\frac{3}{4}}\,dy_N\,ds
\right.\\
& \qquad\qquad\qquad\qquad
\left.
+\int_0^t\int_1^\infty\tau_\varepsilon^{-\frac{1}{2}}\exp\left(-\frac{(x_N-y_N)^2}{4\tau_\varepsilon}\right)\,dy_N\,ds
\right\}\\
& 
\le C\varepsilon^{\frac{1}{2}}\|v\|_{X_T}\left\{t^{\frac{1}{2}}
+\varepsilon^{\frac{1}{2}}\int_0^t(t-s)^{-\frac{1}{2}}s^{\frac{1}{4}}\,ds+t\right\}
\le C\varepsilon^{\frac{1}{2}}\|v\|_{X_T}
(t^{\frac{1}{2}}+\varepsilon^{\frac{1}{2}} t^{\frac{3}{4}}+t)
\end{split}
\end{equation*}
for $x\in\overline{{\mathbb R}^N_+}$ and $0<t<T$,
where $\tau_\varepsilon=\varepsilon^{-1}(t-s)$.
Then, taking a sufficiently small $T>0$ if necessary, we obtain 
\begin{equation}
\label{eq:3.9}
\sup_{0<t<T}\|\tilde{D}_\varepsilon[v](t)\|_{L^\infty}
\le\frac{1}{8}\|v\|_{X_T}. 
\end{equation}
On the other hand,
similarly to \eqref{eq:2.20}, 
by \eqref{eq:3.2} 
we see that
\begin{equation*}
\begin{split}
& \left|\partial_{x_N}\tilde{D}_\varepsilon[v](x,t)\right|\\
& \le 
C\int_0^t\int_0^\infty 
\tilde{K}(x_N,y_N,\tau_\varepsilon)
\|F_2[v](\cdot,y_N,s)\|_{L^\infty({\mathbb R}^{N-1})}\,dy_N\,ds\\
& \le 
C\varepsilon^{\frac{1}{2}}\|v\|_{X_T}\int_0^t\int_0^\infty 
\tilde{K}(x_N,y_N,\tau_\varepsilon)
\left(s^{-\frac{1}{2}}+h(y_N,s)\right)\,dy_N\,ds\\
& \le 
C\varepsilon^{\frac{1}{2}}\|v\|_{X_T}\int_0^t\int_0^\infty 
\tilde{K}(x_N,y_N,\tau_\varepsilon)
\left(s^{-\frac{1}{2}}+y_N^{-\frac{3}{4}}s^{\frac{1}{4}}+y_N^{-\frac{1}{2}}\right)\,dy_N\,ds
\end{split}
\end{equation*}
for $x\in\overline{{\mathbb R}^N_+}$ and $0<t<T$,
where $\tilde{K}$ is the function given by \eqref{eq:2.21}.
By \eqref{eq:2.17} we have
\begin{equation*}
\begin{split}
& \left|\partial_{x_N}\tilde{D}_\varepsilon[v](x,t)\right|
\le 
C\varepsilon^{\frac{1}{2}}\|v\|_{X_T}\bigg(\int_0^ts^{-\frac{1}{2}}\tau_\varepsilon^{-\frac{1}{2}}\,ds
+\int_0^ts^{\frac{1}{4}}\tau_\varepsilon^{-\frac{7}{8}}\,ds
+\int_0^t\tau_\varepsilon^{-\frac{3}{4}}\,ds\bigg)
\\
& =
C\varepsilon^{\frac{1}{2}}\|v\|_{X_T}\bigg(\int_0^ts^{-\frac{1}{2}}(\varepsilon^{-1}(t-s))^{-\frac{1}{2}}\,ds
\\
&\qquad\qquad\qquad
+\int_0^ts^{\frac{1}{4}}(\varepsilon^{-1}(t-s))^{-\frac{7}{8}}\,ds
+\int_0^t(\varepsilon^{-1}(t-s))^{-\frac{3}{4}}\,ds\bigg)
\\
& \le 
C\varepsilon^{\frac{1}{2}}\|v\|_{X_T}\bigg(\varepsilon^{\frac{1}{2}}+\varepsilon^{\frac{7}{8}}t^{\frac{3}{8}}+\varepsilon^{\frac{3}{4}}t^{\frac{1}{4}}\bigg)
\\
&
\le
C(\varepsilon^{-1}t)^{-\frac{1}{2}}\|v\|_{X_T}\bigg((\varepsilon t)^{\frac{1}{2}}+(\varepsilon t)^{\frac{7}{8}}+(\varepsilon t)^{\frac{3}{4}}\bigg)
\end{split}
\end{equation*}
for $x\in\overline{{\mathbb R}^N_+}$ and $0<t<T$.
Taking a sufficiently small $T>0$ if necessary, we see that 
\begin{equation}
\label{eq:3.10}
\sup_{0<t<T}\,(\varepsilon^{-1}t)^{\frac{1}{2}}\|\partial_{x_N}\tilde{D}_\varepsilon[v](t)\|_{L^\infty}
\le\frac{1}{8}\|v\|_{X_T}. 
\end{equation}
Therefore, by \eqref{eq:3.9} and \eqref{eq:3.10} 
we have \eqref{eq:3.8}. 
Furthermore, 
by \eqref{eq:3.2} we apply a similar argument as in the proof of Lemma~\ref{Lemma:2.2}
and deduce that
$\tilde{D}_\varepsilon[v]$ and $\partial_{x_N}\tilde{D}_\varepsilon[v]$ are bounded and smooth in $\overline{{\mathbb R}^N_+}\times(\tau,T)$ for any $0<\tau<T$. 
Thus Lemma~\ref{Lemma:3.2} follows.
$\Box$ \vspace{5pt}

Now we are ready to complete the proof of Theorem~\ref{Theorem:1.1}.
\vspace{5pt}
\newline
{\bf Proof of Theorem~\ref{Theorem:1.1}.}
Let 
\begin{equation}
\label{eq:3.11}
m:=16\max\{\|\varphi\|_{L^\infty}, |\varphi_b|_{L^\infty}\}.
\end{equation}
Let $T_*>0$ be as in Lemma~\ref{Lemma:3.2} 
and $v\in X_{T_*}$ with $\|v\|_{X_{T_*}}\le m$. 
Then, by property~($G_2$), Lemmata~\ref{Lemma:2.2} and~\ref{Lemma:3.2} 
we see that $Q_\varepsilon[v]\in X_{T_*}$. 
Since it follows from \eqref{eq:1.9} and \eqref{eq:2.8} that
$$
\|\Phi\|_{L^\infty}\le\|\varphi\|_{L^\infty}+|\varphi_b|_{L^\infty},
$$
by \eqref{eq:2.1}, \eqref{eq:2.2} and \eqref{eq:3.11} we have
\begin{equation}
\label{eq:3.12}
\|S_1(\varepsilon^{-1}t)\Phi\|_{L^\infty}
+(\varepsilon^{-1} t)^{\frac{1}{2}}\|\partial_{x_N}[S_1(\varepsilon^{-1}t)\Phi]\|_{L^\infty}
\le2\|\Phi\|_{L^\infty}\le\frac{m}{4}
\end{equation}
for $0<t<T_*$. 
Furthermore, by \eqref{eq:2.10} and \eqref{eq:2.19}, 
taking a sufficiently small $T_*>0$ if necessary, we see that 
\begin{equation}
\label{eq:3.13}
\|D_\varepsilon[\varphi_b](t)\|_{L^\infty}
+(\varepsilon^{-1} t)^{\frac{1}{2}}\|\partial_{x_N}D_\varepsilon[\varphi_b](t)\|_{L^\infty}
\le CT_*^{\frac{1}{4}}(1+T_*^{\frac{3}{4}})|\varphi_b|_\infty
\le\frac{m}{4}
\end{equation}
for $0<t<T_*$. 
Lemma~\ref{Lemma:3.2} together with \eqref{eq:3.1}, \eqref{eq:3.12} and \eqref{eq:3.13} implies that 
$$
\left\|Q_\varepsilon[v]\right\|_{X_{T_*}}\le
\|S_1(\varepsilon^{-1}\cdot)\varphi\|_{X_{T_*}}+\|D_\varepsilon[\varphi_b]\|_{X_{T_*}}+\|\tilde{D}_\varepsilon[v]\|_{X_{T_*}}
\le\frac{3m}{4}<m.
$$
Similarly, we obtain 
$$
\left\|Q_\varepsilon[v_1]-Q_\varepsilon[v_2]\right\|_{X_{T_*}}
=\|\tilde{D}_\varepsilon[v_1]-\tilde{D}_\varepsilon[v_2]\|_{X_{T_*}}\le\frac{1}{4}\|v_1-v_2\|_{X_{T_*}}
$$
for $v_i\in X_{T_*}$ with $\|v_i\|_{X_{T_*}}\le m$ $(i=1,2)$. 
Then, the contraction mapping theorem ensures that 
there exists a unique $v_\varepsilon\in X_{T_*}$ with $\|v_\varepsilon\|_{X_{T_*}}\le m$ and
\begin{equation}
\label{eq:3.14}
v_\varepsilon=Q_\varepsilon[v_\varepsilon]
=S_1(\varepsilon^{-1}t)\Phi
-D_\varepsilon[\varphi_b](t)
-\tilde{D}_\varepsilon[v_\varepsilon](t)
\quad\mbox{in}\quad X_{T_*}.
\end{equation}
In particular, we see that
\begin{equation}
\label{eq:3.15}
\|v_\varepsilon\|_{X_{T_*}}\le C(\|\varphi\|_{L^\infty}+|\varphi_b|_{L^\infty}).
\end{equation}
Furthermore, by $(G_2)$ and Lemmata~\ref{Lemma:2.2},~\ref{Lemma:3.2}, we see that $v_\varepsilon$ is bounded and smooth in $\overline{{\mathbb R}^N_+}\times(T_1,T_*)$
for any $0<T_1<T_*$. 
As before, set
$$
w_\varepsilon(x,t)=[S_2(t)\varphi_b](x)+\int_0^t[S_2(t-s)\partial_{x_N}v_\varepsilon(s)](x)\, ds
$$
for $x\in\overline{{\mathbb R}^N_+}$ and $t\in(0,{T_*})$. 
By \eqref{eq:2.8} and \eqref{eq:3.11} we obtain
\begin{equation}
\label{eq:3.16}
\begin{split}
&
\|w_\varepsilon(t)\|_{L^\infty}
\le \|S_2(t)\varphi_b\|_{L^\infty}+\int_0^t\|S_2(t-s)\partial_{x_N}v_\varepsilon(s)\|_{L^\infty}\,ds\\
&\quad
  \le |\varphi_b|_{L^\infty}+\int_0^t|\partial_{x_N}v_\varepsilon(s)|_{L^\infty}\,ds
  \\
  &\quad
 \le \frac{m}{16}+\varepsilon^{\frac{1}{2}}m\int_0^ts^{-\frac{1}{2}}\, ds
 \\
  &\quad
 \le C(1+T_*^{\frac{1}{2}})m
  \le C(1+T_*^{\frac{1}{2}})(\|\varphi\|_{L^\infty}+|\varphi_b|_{L^\infty})<\infty
\end{split}
\end{equation}
for all $0<t<T_*$.
Furthermore, by (${\rm P_3}$) 
we apply a similar argument as in Lemma~\ref{Lemma:2.2} 
and see that $w_\varepsilon$ is bounded and smooth in $\overline{{\mathbb R}^N_+}\times(T_1,T_*)$ for any $0<T_1<T_*$. 
Therefore we deduce that 
$(v_\varepsilon,w_\varepsilon)$ is a solution of \eqref{eq:1.11} in ${\mathbb R}^N_+\times(0,T_*)$. 
In addition, by \eqref{eq:3.15} and \eqref{eq:3.16} 
we have assertion~(a) for any $\tau\in(0,T_*)$.
Since $T_*$ is independent of $m$, 
due to the semigroup properties of $S_1(t)$ and $S_2(t)$, 
we see that $(v_\varepsilon,w_\varepsilon)$ is a global-in-time solution of \eqref{eq:1.11} 
and it satisfies assertion~(a) for any $\tau>0$.

Let $(\tilde{v}_\varepsilon,\tilde{w}_\varepsilon)$ be a global-in-time solution of \eqref{eq:1.11} 
satisfying \eqref{eq:1.14}.
Since 
$$
v_\varepsilon-\tilde{v}_\varepsilon
=Q_\varepsilon[v_\varepsilon]-Q_\varepsilon[\tilde{v}_\varepsilon]
=\tilde{D}_\varepsilon[v_\varepsilon-\tilde{v}_\varepsilon]\quad\mbox{in}\quad X_{T_*},
$$
by \eqref{eq:3.8} we have 
$$
\|v_\varepsilon-\tilde{v}_\varepsilon\|_{X_{T_*}}\le\frac{1}{4}\|v_\varepsilon-\tilde{v}_\varepsilon\|_{X_{T_*}}.
$$
This implies that $v_\varepsilon=\tilde{v}_\varepsilon$ in $X_{T_*}$. 
Repeating this argument, we see that $v_\varepsilon=\tilde{v}_\varepsilon$ in $X_T$ for any $T>0$. 
Therefore we deduce that $(v_\varepsilon,w_\varepsilon)$ is a unique global-in-time solution of \eqref{eq:1.11} 
satisfying \eqref{eq:1.14}.

It remains to prove assertions~(b) and (c). 
Let $T'>0$ and $L>0$.
By \eqref{eq:1.14} and \eqref{eq:2.8}
we have
\begin{equation*}
\begin{split}
&
\|w_\varepsilon(t)-S_2(t)\varphi_b\|_{L^\infty}
\le
\int_0^t\|S_2(t-s)\partial_{x_N}v_\varepsilon(s)\|_{L^\infty}\, ds
\\
&\quad
\le
\int_0^t|\partial_{x_N}v_\varepsilon(s)|_{L^\infty}\,ds
\le C\|v_\varepsilon\|_{X_{T'}}\varepsilon^{\frac{1}{2}}\int_0^ts^{-\frac{1}{2}}\, ds
\le C\|v_\varepsilon\|_{X_{T'}}\varepsilon^{\frac{1}{2}}T'^{\frac{1}{2}}
\end{split}
\end{equation*}
for all $t\in(0,T')$.
This implies assertion~(c).
On the other hand, 
since $\tilde{D}_\varepsilon[v_\varepsilon]$ is given with $F_1[\psi]$ replaced by $F_2[v_\varepsilon]$,
by \eqref{eq:3.2} we apply a similar argument as in \eqref{eq:2.16} 
to obtain
\begin{equation*}
\begin{split}
&
\left|\tilde{D}_\varepsilon[v_\varepsilon](x,t)\right|
\le
C
 \int_0^t\int_0^\infty\bigg(\Gamma_1(x_N-y_N,\tau_\varepsilon)-\Gamma_1(x_N+y_N,\tau_\varepsilon)\bigg)\\
 & \qquad\qquad\qquad\qquad\qquad\qquad\qquad\qquad
 \times
 \|F_2[v_\varepsilon](\cdot,y_N,s)\|_{L^\infty({\mathbb R}^{N-1})}\,dy_N\,ds\\
&
\le C\varepsilon^{\frac{1}{2}}\|v_\varepsilon\|_{X_{T'}}
\int_0^t\int_0^\infty
\bigg(\Gamma_1(x_N-y_N,\tau_\varepsilon)-\Gamma_1(x_N+y_N,\tau_\varepsilon)\bigg)\\
& \qquad\qquad\qquad\qquad\qquad\qquad\qquad\qquad\qquad
 \times
\left(s^{-\frac{1}{2}}+h(y_N,s)\right)\,dy_N\,ds
\\
&
\le C\varepsilon^{\frac{1}{2}}\|v_\varepsilon\|_{X_{T'}}
\int_0^t\int_0^\infty
\bigg(\Gamma_1(x_N-y_N,\tau_\varepsilon)-\Gamma_1(x_N+y_N,\tau_\varepsilon)\bigg)\\
& \qquad\qquad\qquad\qquad\qquad\qquad
 \times
\left(s^{-\frac{1}{2}}+\chi_{y_n\ge1}+y_N^{-\frac{3}{4}}s^{\frac{1}{4}}\chi_{0\le y_N<1}\right)\,dy_N\,ds
\end{split}
\end{equation*}
for $x \in {\mathbb R}^{N-1}\times(0,L)$, $0<t<T'$ and $\varepsilon>0$,
where $\tau_\varepsilon=\varepsilon^{-1}(t-s)$.
By \eqref{eq:2.6} we have
\begin{equation*}
\begin{split}
&
\left|\tilde{D}_\varepsilon[v_\varepsilon](x,t)\right|
\\
&
\le C\varepsilon^{\frac{1}{2}}\|v_\varepsilon\|_{X_{T'}}\left\{\int_0^t(1+s^{-\frac{1}{2}})\tau_\varepsilon^{-\frac{1}{2}}\,ds
+\int_0^t\tau_\varepsilon^{-\frac{1}{2}}s^{\frac{1}{4}}\int_0^1y_N^{-\frac{3}{4}}\,dy_N\,ds
\right\}
\\
&
\le 
C\varepsilon^{\frac{1}{2}}\|v_\varepsilon\|_{X_{T'}}\left\{\int_0^t(1+s^{-\frac{1}{2}})\tau_\varepsilon^{-\frac{1}{2}}\,ds
+\int_0^t\tau_\varepsilon^{-\frac{1}{2}}s^{\frac{1}{4}}\,ds
\right\}
\\
&
\le 
C\varepsilon\|v_\varepsilon\|_{X_{T'}}\left\{\int_0^t(1+s^{-\frac{1}{2}})(t-s)^{-\frac{1}{2}}\,ds
+\int_0^t(t-s)^{-\frac{1}{2}}s^{\frac{1}{4}}\,ds
\right\}
\\
&
\le
C\varepsilon\|v_\varepsilon\|_{X_{T'}}(1+t^{\frac{1}{2}}+t^{\frac{3}{4}})
\end{split}
\end{equation*}
for $x \in {\mathbb R}^{N-1}\times(0,L)$, $0<t<T'$ and $\varepsilon>0$.
This implies that 
\begin{equation}
\label{eq:3.17}
\lim_{\varepsilon\to0}\sup_{t\in(0,T')}\,\|\tilde{D}_\varepsilon[v_\varepsilon](t)\|_{L^\infty({\mathbb R}^{N-1}\times(0,L))}=0.
\end{equation}
Therefore, 
applying \eqref{eq:2.3}, \eqref{eq:2.11} and \eqref{eq:3.17} to \eqref{eq:3.14}, we obtain assertion~(b).
Thus the proof of Theorem~\ref{Theorem:1.1} is complete.
$\Box$\vspace{5pt}

\noindent
{\bf Proof of Corollary~\ref{Corollary:1.1}.}
Corollary~\ref{Corollary:1.1} immediately follows from Theorem~\ref{Theorem:1.1} and Definition~\ref{Definition:1.1}.
$\Box$
\vspace{8pt}

\noindent
{\bf Acknowledgment.} The first author was supported in part by the Slovak
Research and Development Agency under the contract No. APVV-14-0378 and by the VEGA grant
1/0347/18. Part of this work was carried out while the first author visited
the Research Alliance Center for Mathematical Sciences, Tohoku University. 
The second author was supported in part by the Grant-in-Aid for Scientific Research (A)(No.~15H02058)
from Japan Society for the Promotion of Science.
The third author was supported by the Grant-in-Aid for Young Scientists (B)
(No.~16K17629)
from Japan Society for the Promotion of Science. 
%


\end{document}